\documentclass[a4paper]{gtart}
\usepackage{amsmath, amssymb, verbatim}
\usepackage{pinlabel, subfig}
\usepackage{graphs}
\usepackage{graphicx}
\newtheorem{thm}{Theorem}[section]
\newtheorem{cor}[thm]{Corollary}
\newtheorem{lem}[thm]{Lemma}
\newtheorem{exa}[thm]{Example}
\newtheorem{prop}[thm]{Proposition}

\theoremstyle{definition}
\newtheorem{defn}[thm]{Definition}

\theoremstyle{remark}
\newtheorem{rem}[thm]{Remark}
\newtheorem{rems}[thm]{Remarks}
\numberwithin{equation}{section}

\newcommand{\al}{\alpha}
\newcommand{\be}{\beta}
\newcommand{\de}{\delta}

\newcommand{\si}{\sigma}

\newcommand{\la}{\lambda}
\newcommand{\Si}{\Sigma}

\renewcommand{\th}{\theta}

\newcommand{\ka}{\kappa}
\newcommand{\ta}{\tau}

\renewcommand{\t}{\mathbf t}

\newcommand{\s}{\mathbf s}
\newcommand{\n}{\mathbf n}

\newcommand{\LL}{\mathcal L}
\newcommand{\FF}{\mathcal F}

\newcommand{\Z}{\mathbb Z}
\newcommand{\Q}{\mathbb Q}
\newcommand{\R}{\mathbb R}
\newcommand{\T}{\mathbb T}

\newcommand{\CP}{{\mathbb C}{\mathbb P}}
\newcommand{\del}{\partial}

\newcommand{\x}{\times}
\newcommand{\lra}{\longrightarrow}

\newcommand{\wHF}{\widehat{HF}}
\newcommand{\Li}{\mathbb {L}}
\DeclareMathOperator{\tb}{tb}
\DeclareMathOperator{\sel}{sl}

\DeclareMathOperator{\rot}{rot}

\DeclareMathOperator{\Cont}{Cont}
\newcommand{\Spin}{\mathrm{Spin}}

\begin{document}

\title{Contact surgery and transverse invariants} 

\author{Paolo Lisca}
\address{Dipartimento di Matematica ``L. Tonelli''\\ 
Largo Bruno Pontecorvo, 5\\
Universit\`a di Pisa \\
56127 Pisa, ITALY} 
\email{lisca@dm.unipi.it}

\author{Andr\'{a}s I. Stipsicz}
\address{R\'enyi Institute of Mathematics\\
Hungarian Academy of Sciences\\
H-1053 Budapest\\ 
Re\'altanoda utca 13--15, Hungary}
\email{stipsicz@math-inst.hu}

\date{}

\begin{abstract} 
The purpose of this paper is two--fold: (1) to derive new existence results for tight 
contact structures on closed 3--manifolds presented by integral surgery along 
knots in $S^3$, and (2) to introduce a new invariant for transverse knots 
in contact 3--manifolds. Regarding (1), we extend our previous existence 
results from surgeries along knots of genus $g$ and maximal 
Thurston--Bennequin number $2g-1$ to surgeries along knots of 
genus $g$ and maximal self--linking number $2g-1$.
\end{abstract}

\maketitle

\section{Introduction}\label{s:intro}
One of the motivating questions of 3--dimensional contact topology is
to characterize those (closed, oriented) 3--manifolds which admit
(positive, cooriented) tight contact structures. This question was
answered recently for Seifert fibered 3--manifolds~\cite{Duke}, but
the general case is still wide open. A particular family of
3--manifolds is given by those which can be presented as surgery along
a knot in $S^3$. It seems reasonable to expect that the use of the
invariant $\widehat\LL(\ka)$ for Legendrian knots in $S^3$ defined
in~\cite{LOSS} together with contact geometric constructions might
provide a way to find tight examples on many such surgeries. Some
justification for such an expectation is provided by a result of
Sahamie \cite{bijan}, showing that if the Legendrian invariant
$\widehat\LL(\ka)$ of a Legendrian knot $\ka$ vanishes, then the contact
Ozsv\'ath--Szab\'o invariant $c(\xi_1(\ka))$ of the result of contact
$(+1)$--surgery along $\ka$ is also zero. On the other hand, for a
Legendrian knot $\ka$ in the standard contact $S^3$ satisfying
$\tb(\ka )=2g_s(\ka )-1>0$ (where $g_s(\ka )$ denotes the smooth
4--ball genus of the knot type of $\ka$) it was shown in \cite{LS1}
that the result of contact $(+1)$--surgery has nonvanishing contact
Ozsv\'ath--Szab\'o invariant, implying in particular tightness for the
contact structure. (The nonvanishing of this invariant implies
tightness, while a contact structure with vanishing invariant might be
either tight or overtwisted.)

In this paper we extend this nonvanishing result to knots with other 
properties, allowing contact surgeries with higher coefficients.
Given a knot type $K\subset S^3$, let the \emph{maximal self--linking number}
of $K$ be the largest self--linking number of a transverse representative 
of $K$ (with respect to the standard contact structure $\xi_{st}$). 
Also, denote the Seifert genus of $K$ by $g(K)$. 

\begin{thm}\label{t:mainsl}
Let $K\subset S^3$ be a knot type with maximal self--linking number equal to 
$2g(K)-1$. Then, for $r\geq 2g(K)$ the 3--manifold $S^3_{r}(K)$ carries 
tight contact structures.
\end{thm}

Examples of knots satisfying the assumptions of the above theorem are
provided by strongly quasi-positive, fibered knots in $S^3$.  In
particular, iterated torus knots $K((p_1,q_1), \ldots, (p_k, q_k))$
with all $p_i, q_i>0$ are such examples.

According to \cite{EH2} the (2,3)--cable $K_{2,3}$ of the (2,3) torus
knot $T_{2,3}$ provides an example of a knot for which
Theorem~\ref{t:mainsl} applies while the previous result from
\cite{LS1} does not: the maximal self--linking of $K_{2,3}$ is equal
to 7 (which is equal to $2g(K_{2,3}) -1$) while the maximal
Thurston--Bennequin number of $K_{2,3}$ is 6. By taking the connected
sum of $n$ copies of this knot, the difference between the maximal
self--linking and the maximal Thurston--Bennequin number can be made,
in fact, arbitrarily large. Related, prime knot examples for the same
phenomenon are provided by $(p,q)$--cables ($q>p\geq 1$) $K_{p,q}$ of
the (2,3) torus knot $T_{2,3}$: according to \cite{bulent} the maximal
self-linking number of $K_{p,q}$ (which again coincides with
$2g(K_{p,q})-1$) is equal to $pq+q-p$, while the maximal
Thurston--Bennequin number of $K_{p,q}$ is $pq$.

\bigskip

We found it convenient to organize the surgery theoretic information
about a Legendrian (and about a transverse) knot into an
invariant which takes its values in Heegaard Floer homology
groups (and ultimately in the inverse limit of some of these groups).
Although the resulting surgery invariant $\tilde c$ shares a number of
properties with the Legendrian (and transverse) knot invariants
introduced in \cite{LOSS}, we found a vanishing result for $\tilde c$
(given in Theorem~\ref{t:main}) which is, according to a recent result
of Vela-Vick \cite{SVV}, in sharp contrast with the corresponding
behaviour of the Legendrian invariant $\widehat \LL$ of \cite{LOSS}.

In order to state our results we need some preliminary notation.
Let $Y$ be a closed, oriented 3--manifold and $K\subset Y$ a knot
type.  Let $\FF_K$ be the set of framed isotopy classes of framed
knots in the (unframed) knot type $K$. We will follow the usual practice of 
referring to the elements of $\FF_K$ as to the ``framings'' of $K$.  
Recall that for $K$ null--homologous $\FF_K$ is an affine 
$\Z$--space, and that even if $K$ is not null--homologous this is still true 
if $Y$ is not of the form $Y'\# S^1\x S^2$~\cite{Ch, HP}.
For $k\in\Z$ and $f\in\FF_K$, we shall denote the result of $k$ acting on $f$ by
$f+k$. When $\FF_K$ is an affine $\Z$--space, 
$\FF_K$ inherits a natural linear 
order from $\Z$: if $f, g\in\FF_K$ with $f=g+k$, $k\in\Z$, then $f\geq g$ 
if and only if $k\geq 0$. We will denote by $Y_f(K)$ the 3--manifold resulting 
from surgery on $Y$ along $K$ with framing $f$.

Given a contact 3--manifold $(Y,\xi)$, a {\em framed Legendrian knot}
in $(Y,\xi)$ is a pair $(\ka,f)$, where $\ka\subset (Y,\xi)$ is a
Legendrian knot and $f\in\FF_K$ is a framing of the topological type
$K$ of $\ka$. A {\em framed transverse knot} in $(Y,\xi)$ is a pair
$(\tau,f)$, where $\tau\subset (Y,\xi)$ is a transverse knot and
$f\in\FF_K$ is a framing of the topological type $K$ of $\tau$. Denote
by $\T(Y,\xi,K,f)$ the set of transverse isotopy classes of framed
transverse knots $(\tau,f)$ in $(Y,\xi)$ with $\tau$ in the
topological type $K$. Let $\Cont(Y)$ be the set of isomorphism
classes of contact structures on $Y$.
Fix a transverse knot $\ta \subset (Y, \xi )$ in the knot type
$K$. By considering  a Legendrian approximation $\ka$ of $\ta$, and by
applying appropriate contact surgery along $\ka$ (where the exact meaning
of 'appropriate' will be clarified in Subsection~\ref{ss:cosu}), a contact structure
$\tilde I(\xi , \ta , f)$ can be defined on the 3--manifold
$Y_f(K)$.
\begin{thm}\label{t:maintransverse}
  Let $Y$ be a closed, oriented 3--manifold and $K$ a knot type in
  $Y$.  Suppose that either $K$ is null--homologous or $Y$ is not of
  the form $Y'\# S^1\x S^2$. Given a contact structure $\xi$ on $Y$
  and a framing $f$ on $K\subset Y$, there is a well--defined map
\[
\begin{matrix}
\T(Y,\xi,K,f) & \lra & \Cont(Y_f(K))\\
[(\ta, f)] & \longmapsto & \tilde I(\xi,\ta,f)
\end{matrix}
\]
\end{thm} 
In ~\cite{OSz-cont} Ozsv\'ath and Szab\'o associated an element of the
Heegaard Floer group $\wHF(-Y)$ to every contact 3--manifold $(Y,\xi)$.  By
fixing an identification between the diffeomorphic 3--manifolds $Y_f(\ta)$ and
$Y_f(K)$ we get a family of Heegaard Floer elements
\begin{equation*}
\tilde c(\xi,\ta,f):= c(\tilde I(\xi,\ta,f))\in \wHF(-Y_f(K)),\quad  
f\in\FF_K
\end{equation*}
for every transverse knot $\ta \subset (Y, \xi )$ representing the
knot type $K$. (In this paper we always consider Heegaard Floer
homology with $\Z/2\Z$ coefficients). The elements themselves might
depend on the chosen identification of $Y_f(\ta)$ with $Y_f(K)$; their
vanishing/nonvanishing properties, on the other hand, are independent
of this choice. Since in the following we will exclusively focus on
vanishing/nonvanishing questions, we shall not mention the above
identification again.

The invariant $\tilde c$ is non--trivial. In fact,  
in Example~\ref{exa:ketto} we show, using the main result
of~\cite{LS1}, that if $\ta$ is the link of an isolated curve
singularity in the standard contact 3--sphere $(S^3,\xi_{\rm st})$
then $\tilde c(\xi_{\rm st},\ta, f_S+2g_s(K))\neq 0$,
where $f_S$ is the framing defined by a Seifert surface of
$K$, and $g_s(K)$ is the slice genus of $K$. In Section~\ref{s:def} we
also show (Corollary~\ref{c:defin}) that if $\tilde c(\xi,\ta,f)\neq 0$ then
$\tilde c(\xi,\ta,g)\neq 0$ for every $g\geq f$.

By using appropriate cobordisms and maps induced by them, an inverse
limit $H(Y,K)$ of Heegaard Floer groups of results of surgeries of $Y$ along
$K$ can be defined, and we show that the family 
$\left(\tilde c (\xi , \ta , f)\right)_{f\in\FF_K}$ defines a
single element $\tilde c (\xi , \ta )$ in this limit 
(see Proposition~\ref{p:invlimit}). Notice that
$\tilde c$ is defined for a transverse knot $\ta$ through its
Legendrian approximations, a feature similar to the definition of the
transverse invariant $\widehat {\mathcal {T}}$ of \cite{LOSS} (resting on
the corresponding Legendrian invariant $\widehat \LL$).  For $\tilde c$,
however, we have the following vanishing result which shows, in
particular, that $\tilde c(\xi,\ta )$ behaves quite differently from the
transverse invariant $\widehat {\mathcal {T}}$ of~\cite{LOSS}.

\begin{thm}\label{t:main}
  Let $\Si$ be an oriented surface with boundary and
  $\phi\co\Si\to\Si$ an orientation--preserving diffeomorphism which
  restricts to the identity on a collar around $\del\Si$. Let
  $(Y,\xi_{(\Si,\phi)})$ be the contact 3--manifold compatible with
  the open book decomposition induced by $(\Si,\phi)$. Suppose that
  $b_1(Y)=0$ and $c(\xi_{(\Si,\phi)}) = 0$, and let 
$\tau\subset (Y,\xi_{(\Si,\phi)})$ be a
component of the boundary of $\Si$ viewed as the binding of the open
book. Then, $\tilde c(\xi_{(\Si,\phi)},\tau)=0$.
\end{thm}

Theorem~\ref{t:main} should be contrasted with the main result
of~\cite{SVV}, which says that the transverse invariant $\widehat
{\mathcal {T}}$ of~\cite{LOSS} is nonvanishing for a binding of an open
book. (See also \cite{ESVV} for the case of disconnected bindings.)

The following nonvanishing result provides the desired construction of
tight contact structures on certain surgered 3--manifolds, and will
serve as the main ingredient in the proof of Theorem~\ref{t:mainsl}.

\begin{thm}\label{t:main2}
Suppose that the open book decomposition induced by $(\Si,\phi)$ is
compatible with the standard contact structure $\xi_{st}$ on $S^3$.
Let $\ta\subset S^3$ be a binding component of the open book
decomposition having knot type $K$ and self--linking number
$\sel(\tau)=2g(K)-1$, where $g(K)$ is the Seifert genus of $K$. Then,
$\tilde c(\xi_{st},\ta )\neq 0$. In fact, if $f_S$ denotes the Seifert
framing of $\ta$, the Heegaard Floer homology element $\tilde
c(\xi_{st},\ta, f)$ is nonzero for each $f\geq f_S+2g(K)$.
\end{thm}

The paper is organized as follows. In Section~\ref{s:contact} we
establish the properties of contact surgeries that we use to define
the transverse invariants. In Section~\ref{s:def} we define the
invariants, thus establishing Theorem~\ref{t:maintransverse}, and we
prove their basic properties.  In Section~\ref{s:main} we prove the
vanishing Theorem~\ref{t:main} while in Section~\ref{s:main2} we give
the proofs of the nonvanishing result given by Theorem~\ref{t:main2}
and ultimately we prove Theorem~\ref{t:mainsl}.

{\bf Acknowledgements:} We would like to thank John Etnyre and 
Matt Hedden for stimulating discussions, and the anonymous referee 
for useful suggestions which helped to improve the presentation. 
Part of this work was carried out while the authors visited the Mathematical Sciences Research Institute, Berkeley, as participants of the `Homology theories for knots and links'
special semester. The present work is part of the authors' activities within CAST, 
a Research Network Program of the European Science Foundation. PL was partially 
supported by PRIN 2007, MIUR. AS was partially supported by OTKA Grant NK81203
and by the \emph{Lend\"ulet program}. 

\section{Contact surgeries and stabilizations}\label{s:contact} 
\label{sec:second}

In this section we establish the properties of certain contact surgeries 
which will allow us to define the invariants and study their basic properties. 

Let $K\subset Y$ be a knot type in the closed 3--manifold $Y$.  Let
$\xi$ be a contact structure on $Y$ and $\ka\subset(Y,\xi)$ a
Legendrian knot belonging to $K$. Recall that, given a non--zero
rational number $r\in\Q$, one can perform~\emph{contact $r$--surgery}
along $\ka$ to obtain a new contact 3--manifold $(Y',\xi')$~\cite{DG}.
When $r=\pm 1$ the contact structure $\xi'$ is uniquely determined,
therefore in this case we can safely use the notation $\xi_{\pm  1}(\ka)$ 
for $\xi'$. In general, there are several possible choices
for $\xi'$.  According to~\cite[Proposition~7]{DG}, for
$\frac{p}{q}>1$ every contact $\frac{p}{q}$--surgery on $\ka$ is
equivalent to a contact $(+1)$--surgery on $\ka$ followed by a contact
$\frac{p}{q-p}$--surgery on a Legendrian pushoff copy of $\ka$.
Moreover, by~\cite[Proposition~3]{DG} (see~\cite{DGS} as well) every
contact $r$--surgery along $\ka\subset (Y,\xi)$ with $r<0$ is
equivalent to a Legendrian (i.e. $(-1)$--) surgery along a Legendrian link $\Li
=\cup_{i=0}^m L_i$ belonging to a set determined via a simple
algorithm by the Legendrian knot $\ka$ and the contact surgery
coefficient $r$.  The algorithm to obtain the set of $\Li$'s is the
following. Let
\[
[a_0,\ldots,a_m]:= 
a_0 - \cfrac{1}{a_1 -
       \cfrac{1}{\ddots -
        \cfrac{1}{a_m}
}},
\quad a_0,\ldots a_m\geq 2,
\]
be the continued fraction expansion of $1-r$. To obtain $L_0$, stabilize $a_0-2$ times 
a Legendrian push--off of $\ka$  in every possible way. 
Then, stabilize $a_1-2$ times a Legendrian push--off of $L_0$ in every possible way. 
Repeat the above scheme for each of the remaining pivots of the continued fraction expansion. 

\begin{figure}[h!]
\labellist
\small\hair 2pt
\pinlabel $x$ at 13 140
\pinlabel $\ka$ at 74 114
\pinlabel $\theta$ at 109 90
\pinlabel $x$ at 235 179
\pinlabel $\ka_-$ at 330 166
\pinlabel $\theta$ at 360 130
\pinlabel $\theta$ at 352 17
\pinlabel $x$ at 234 70
\pinlabel $\ka_+$ at 283 58
\endlabellist
\centering
\includegraphics[scale=0.5]{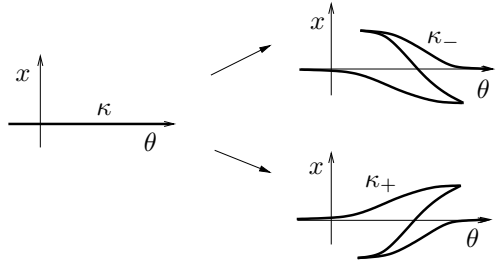}
\caption{Negative and positive Legendrian stabilizations}
\label{f:stab}
\end{figure}

We are interested in contact $n$--surgeries, where $n$ is a positive integer. 
In this case, since 
\[
1-\frac{n}{1-n}=\frac{2n-1}{n-1}=[3,\overbrace{2,\ldots,2}^{n-2}], 
\]
there are only two choices for the stabilizations of $\ka$, because the choice of the first 
one determines all the others. An orientation of $\ka$ allows one to specify 
unambiguously such a choice, because it specifies the~\emph{negative stabilization} 
$\ka_-$ and the~\emph{positive stabilization} $\ka_+$ of $\ka$. In
a standard neighborhood $\R/\Z\x \R^2$ of $\ka$ with coordinates $(\th,x,y)$ the 
contact structure is given by $\xi=\ker(dx+yd\th)$, and the $(\th,x)$--projections of 
$\ka$, $\ka_-$ and $\ka_+$ are illustrated in Figure~\ref{f:stab}. 
{From} now on we shall assume that every Legendrian knot $\ka$ is oriented. 
\begin{defn}\label{d:xi}
We denote by $\xi_n^-(\ka)$ (respectively $\xi_n^+(\ka)$)   
the contact structure corresponding to the choice of the~\emph{negative} 
stabilization $\ka_-$ (respectively the~\emph{positive} stabilization $\ka_+$) of 
the oriented Legendrian knot $\ka$ (see Figure~\ref{f:stab}). 
\end{defn} 

Observe that $\ka_+$ and $\ka_-$ inherit an orientation from $\ka$ in a natural way. 
We shall always assume that $\ka_+$ and $\ka_-$ are given the  orientation 
induced by $\ka$. 

\begin{lem}\label{l:changesign}
Let $\ka\subset (Y,\xi)$ be an oriented Legendrian knot. Then, 
for each $n>0$ we have 
\[
\xi_n^+(\ka) = \xi_n^-(-\ka).
\]
\end{lem} 

\begin{proof} 
The statement follows from the definition of contact surgery 
together with the easily checked fact that $(-\ka)_-=-(\ka_+)$ for every 
oriented Legendrian knot $\ka$. 
\end{proof} 

We want to study the contact structure $\xi^-_n(\ka)$ when 
$\ka$ is a stabilization. The following lemma was proved in greater generality in~\cite{Ozb} using 
the main result of~\cite{HKM}. Here we give  a simple and constructive proof. 

\begin{lem}\label{l:overtwisted}
Let $\ka\subset (Y,\xi)$ be an oriented Legendrian knot. Then, 
$\xi_n^-(\ka_+)$ is an overtwisted contact structure for each $n>0$.
\end{lem}

\begin{proof}
Ozbagci~\cite[Proposition~13]{Ozb} shows that for $r>0$ any contact 
$r$--surgery on a positive stabilization in which the Legendrian 
pushoffs are all negative stabilizations is overtwisted, by 
constructing a non--right veering compatible open book 
and appealing to the results of~\cite{HKM}. 

On the other hand, the lemma can be easily checked directly as follows. 
The left--hand side of Figure~\ref{f:otdisk} illustrates the 
contact surgery yielding $\xi_n^-(\ka_+)$ in a standard neighborhood 
of $\ka$.
\begin{figure}[h!]
\labellist
\small\hair 2pt
\pinlabel {\scriptsize $+1$} at 19 63
\pinlabel $\ka_+$ at 171 62
\pinlabel {\scriptsize $-1$} at 209 88
\pinlabel {\scriptsize $-1$} at 223 114
\pinlabel $\ka$ at 243 17
\pinlabel {\scriptsize $+1$} at 451 81
\pinlabel {\scriptsize $-1$} at 519 93
\pinlabel {\scriptsize $-1$} at 536 120
\pinlabel $\ka_+$ at 573 62
\pinlabel $\tilde\ka$ at 573 17
\endlabellist
\centering
\includegraphics[scale=0.5]{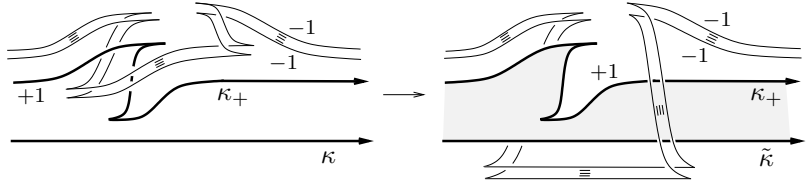}
\caption{The overtwisted disk in $\xi_n^-(\ka_+)$. Notice that the knots labeled by $\ka$
on the left and $\tilde\ka$ on the right are not necessarily isotopic in the surgered manifold. 
Nevertheless, $\ka _+$ and $\tilde\ka$ provide the shaded annulus, which 
then caps off to an overtwisted disk with boundary equal to $\tilde\ka$.}
\label{f:otdisk}
\end{figure}
The right--hand side of the picture shows how the $n-1$ push--offs 
of $(\ka_+)_-$ can be Legendrian isotoped until one can see 
the shaded overtwisted disk. 
\end{proof} 

The following proposition gives the key property of $\xi_n^-(\ka)$ which 
yields transverse invariants. 

\begin{prop}\label{p:stab}
  Let $\ka _-$ denote the negative Legendrian stabilization of the
  oriented Legendrian knot $\ka \subset (Y,\xi)$. Then, for each $n>0$
  the contact structures $\xi_{n+1}^-(\ka_-)$ and 
  $\xi_n^-(\ka)$ are isomorphic.
\end{prop}

\begin{figure}[h!]
\labellist
\footnotesize\hair 2pt
\pinlabel $1$ at 114 184
\pinlabel $1$ at 369 185
\pinlabel $2$ at 80 122
\pinlabel $2$ at 338 122
\pinlabel $3$ at 78 70
\pinlabel $3$ at 338 70
\pinlabel $4$ at 112 3
\pinlabel $4$ at 388 8
\pinlabel $+$ at 14 147
\pinlabel $+$ at 49 131
\pinlabel $+$ at 374 135
\pinlabel $-$ at 15 96
\pinlabel $+$ at 53 81
\pinlabel $+$ at 366 93
\pinlabel $-$ at 334 16
\endlabellist
\centering
\includegraphics[scale=0.7]{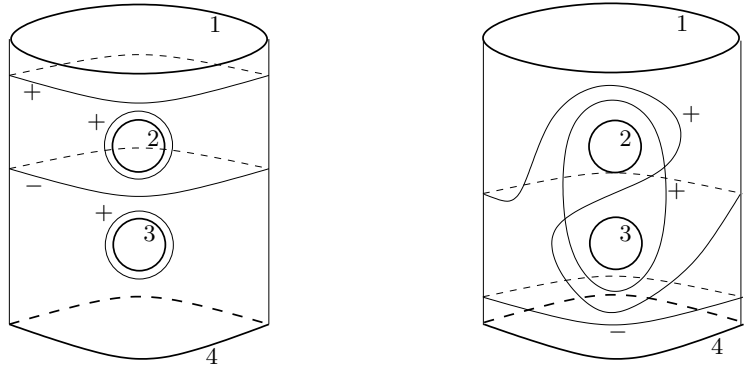}
\caption{The relation $\de_{12}^{-1}\de_1\de_2\de_3 =
  \de_{13}\de_{23}\de_4^{-1}$. Signs on the curves indicate
  whether right-handed $(+)$ and left-handed $(-)$ Dehn twists are to be performed.}
\label{f:lantern}
\end{figure}

Before proving the proposition we recall the lantern relation. Let $A$
be a surface with boundary homeomorphic to a twice punctured annulus.
If we denote by $\de_i$ the positive Dehn twist along a curve parallel
to the $i$--th boundary component of $A$, and by $\de_{ij}$ the
positive Dehn twist along a curve encircling the boundary components
$i$ and $j$, the \emph{lantern relation} reads
$\de_1\de_2\de_3\de_4=\de_{12}\de_{13}\de_{23}$.  In the proof of
Proposition~\ref{p:stab} we are going to use the equivalent relation
$\de_{12}^{-1}\de_1\de_2\de_3 = \de_{13}\de_{23}\de_4^{-1}$.
Figure~\ref{f:lantern} provides a graphical representation of this
relation. In fact, whenever the twice punctured annulus embeds into a surface,
the Dehn twists corresponding to the images of the curves on the diagram satisfy the
lantern relation.

\begin{proof}[Proof of Proposition~\ref{p:stab}]
Consider an open book for $\xi$ with a page  which contains $\ka$ and such that 
the page framing induced on $\ka$ is equal 
to the contact framing of $\ka$. After two Giroux stabilizations we 
can accomodate $\ka$, $\ka_-$ and $(\ka_-)_-$ on the same page of the 
resulting open book, still with equal page and contact framings (see eg~\cite{Et0}). 
After performing a negative Dehn twist along $\ka_-$ and positive Dehn twists 
along $n$ parallel copies of $(\ka_-)_-$ we obtain an open book for 
$\xi^-_{n+1}(\ka_-)$, as illustrated in Figure~\ref{f:isotopy1}.
\begin{figure}[ht!]
\centering
\subfloat[Open book for $\xi_{n+1}^-(\ka_-)$]{
\labellist
\small\hair 2pt
\pinlabel $n$ at 480 72
\pinlabel $\ka$ at 109 15
\pinlabel $\ka_-$ at 113 67
\pinlabel $(\ka_-)_-$ at 150 137
\pinlabel $+$ at 255 164
\pinlabel $+$ at 294 164
\pinlabel $+$ at 317 192
\pinlabel $+$ at 256 84
\pinlabel $-$ at 309 90
\pinlabel {\Large $\}$} at 465 72
\endlabellist
\includegraphics[scale=0.5]{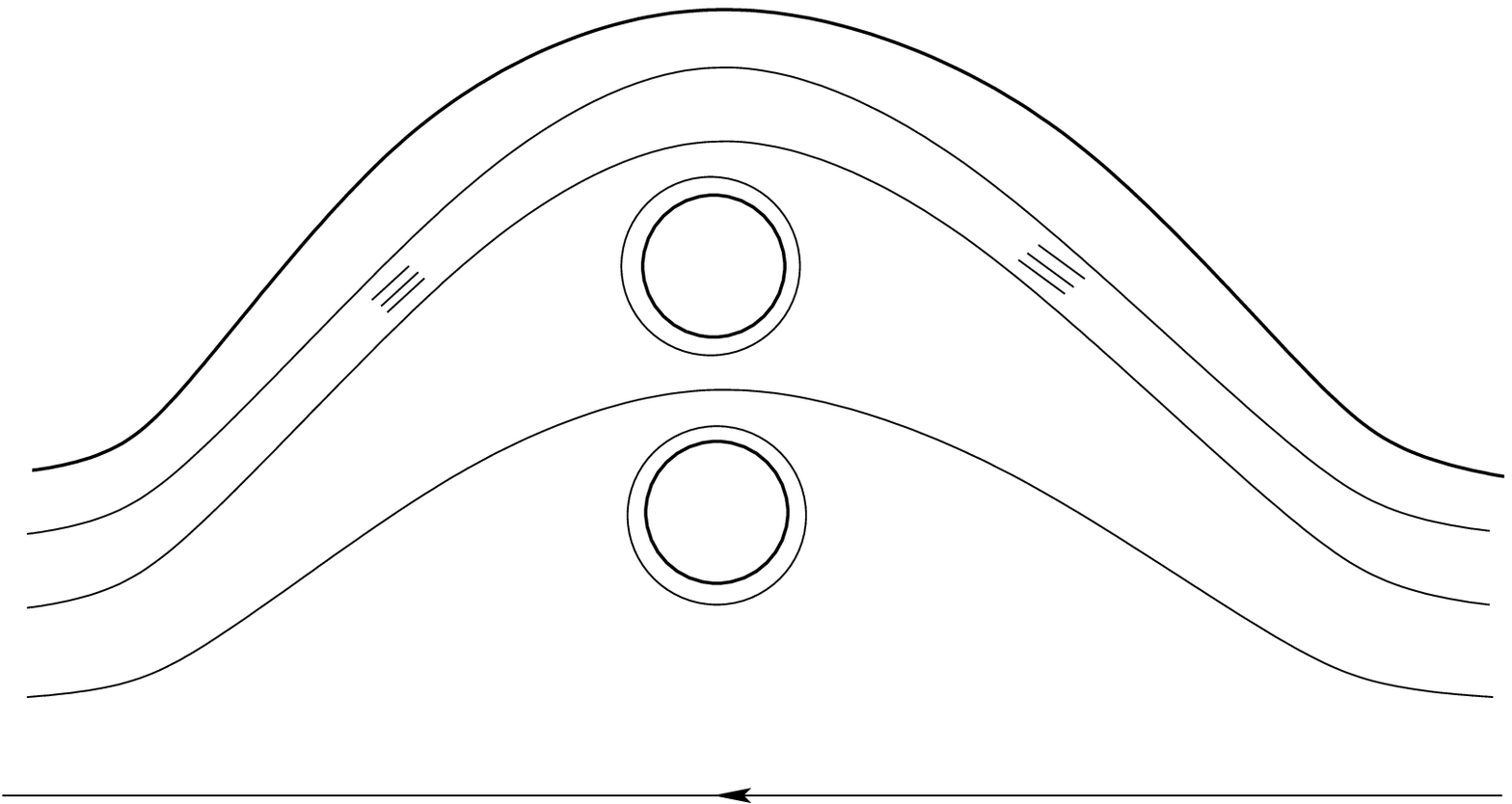}
\label{f:isotopy1}}\\
\subfloat[Open book for $\xi_{n+1}^-(\ka_-)$ after applying the relation]{
\labellist
\small\hair 2pt
\pinlabel $+$ at 314 197
\pinlabel $+$ at 288 165
\pinlabel $+$ at 270 116
\pinlabel $+$ at 315 45
\pinlabel $n-1$ at 500 73
\pinlabel $-$ at 128 12
\pinlabel {\Large $\}$} at 464 73
\endlabellist
\includegraphics[scale=0.45]{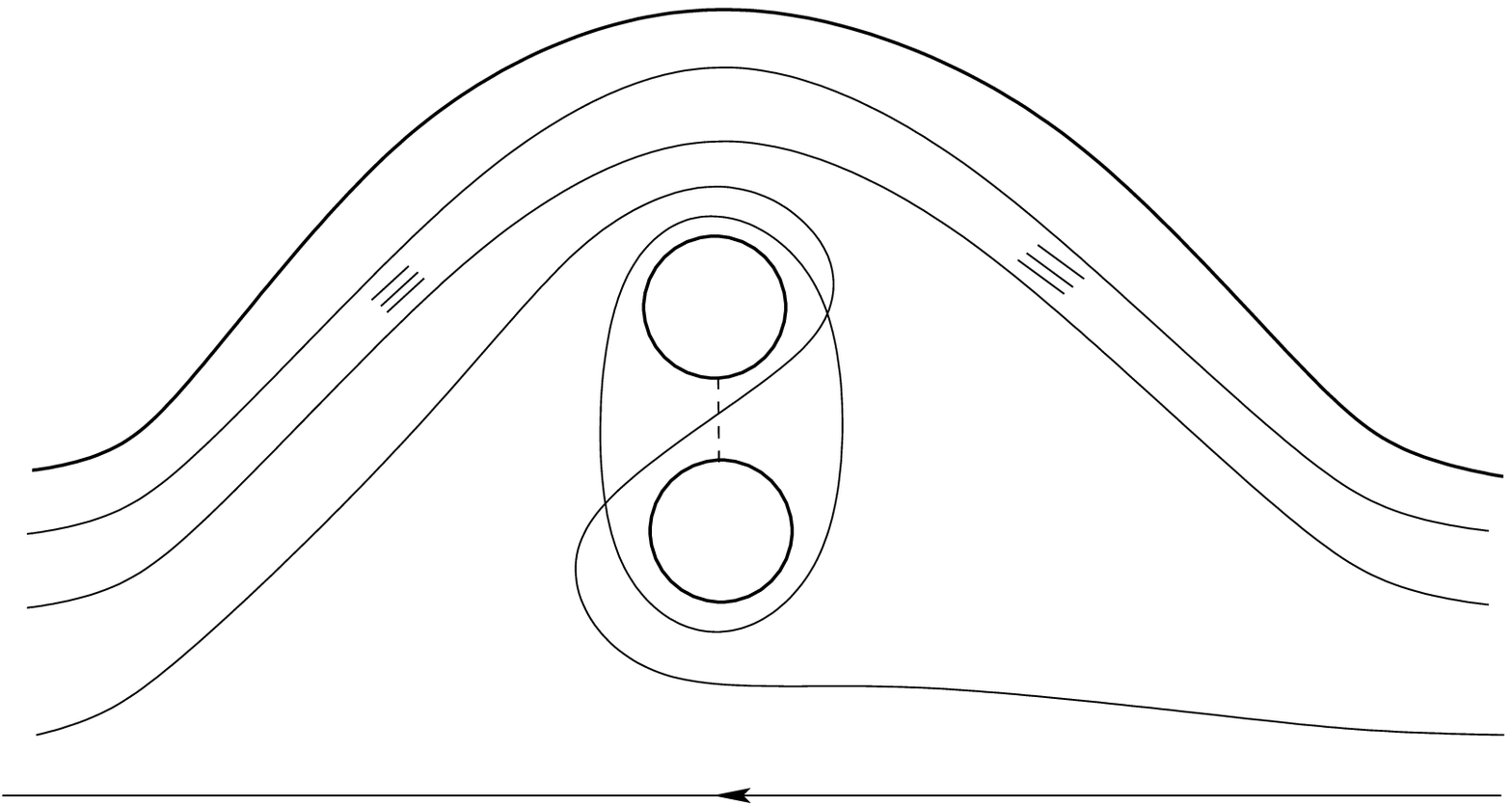}
\label{f:isotopy2}}\\
\subfloat[Open book for $\xi_n^-(\ka)$]{
\labellist
\small\hair 2pt
\pinlabel $+$ at 318 194
\pinlabel $+$ at 289 164
\pinlabel $+$ at 290 105
\pinlabel $n-1$ at 496 72
\pinlabel $-$ at 106 13
\pinlabel $\ka$ at 366 13
\pinlabel $\ka_-$ at 390 69
\pinlabel {\Large $\}$} at 466 72
\endlabellist
\includegraphics[scale=0.5]{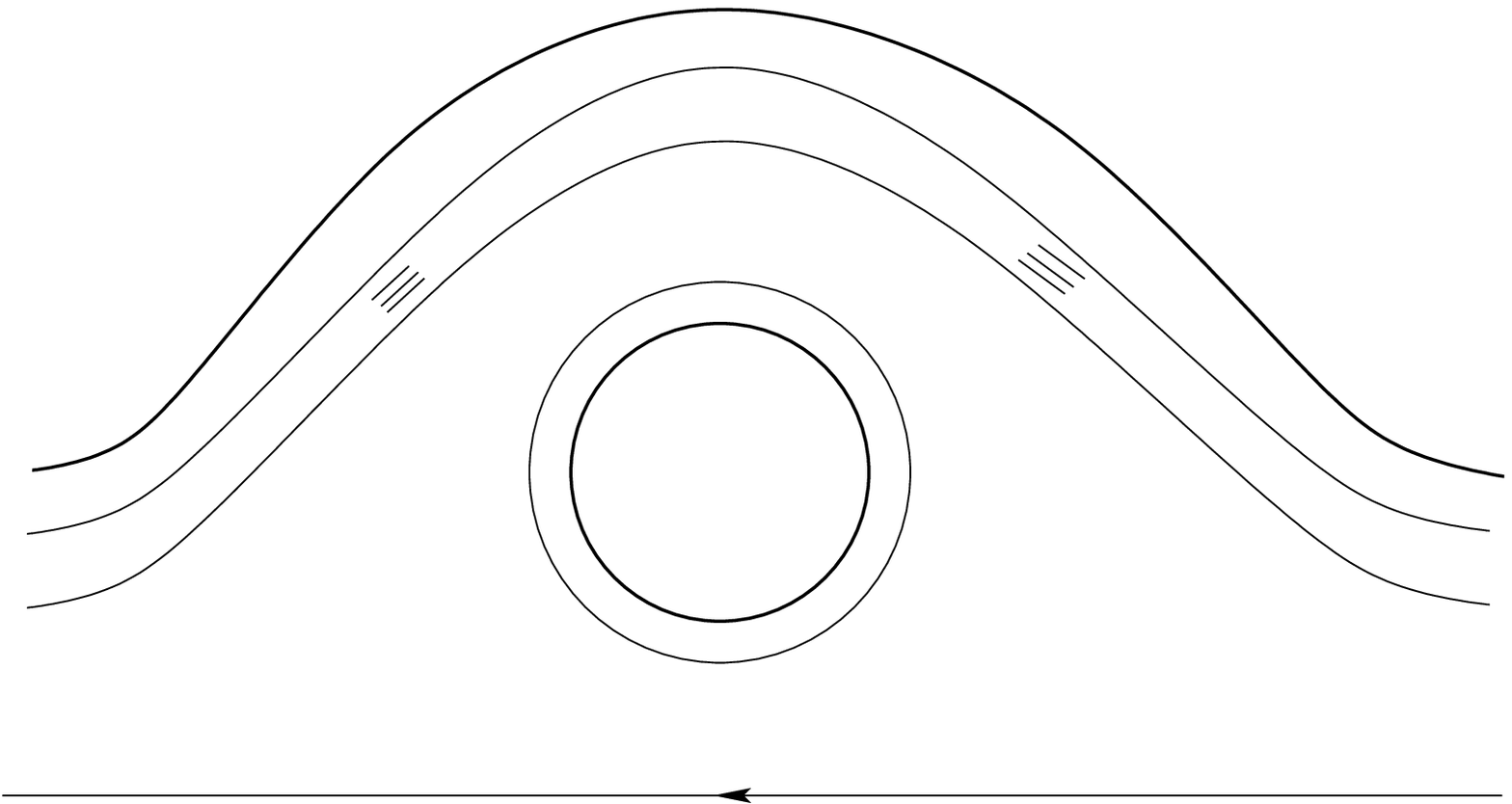}
\label{f:isotopy4}}
\caption{Isomorphism between $\xi_{n+1}^-(\ka_-)$ and $\xi_n^-(\ka)$}
\end{figure}
In Figure~\ref{f:isotopy2} we see what happens to the open book for
$\xi_{n+1}^-(\ka_-)$ when we apply the relation of 
Figure~\ref{f:lantern} inside the
twice punctured annulus visible in the picture.
The dashed arc of Figure~\ref{f:isotopy2} shows that the open book can
be Giroux destabilized, yielding Figure~\ref{f:isotopy4}, which is an
open book for $\xi_n^-(\ka)$.
\end{proof}

The following corollary can be viewed as a generalization
of~\cite[Theorem 1]{Et}.

\begin{cor}\label{c:stab}
  Let $\ka_1, \ka_2\subset (Y,\xi)$ be two Legendrian knots. If after
  negatively stabilizing the same number of times $\ka_1$ and $\ka_2$
  become Legendrian isotopic, then $\xi^-_n(\ka_1)$ is isomorphic to
  $\xi^-_n(\ka_2)$ for each $n>0$.
\end{cor} 

\begin{proof} 
  Suppose that $\ka'_1$ and $\ka'_2$ are Legendrian isotopic
  Legendrian knots obtained by negatively stabilizing $\ka_1$ and
  $\ka_2$ $m$ times.  Then, for each $n>0$, the contact structure
  $\xi^-_{n+m}(\ka'_1)$ is isotopic to $\xi^-_{n+m}(\ka'_2)$.
  Applying Proposition~\ref{p:stab} $m$ times we conclude that
  $\xi^-_{n+m}(\ka'_1)$ is isomorphic to $\xi^-_n(\ka_1)$ and
  $\xi^-_{n+m}(\ka'_2)$ is isomorphic to $\xi^-_n(\ka_2)$.  Therefore
  $\xi^-_n(\ka_1)$ and $\xi^-_n(\ka_2)$ are isomorphic for each $n>0$.
\end{proof} 

Lemma~\ref{l:overtwisted} and Proposition~\ref{p:stab} admit slight refinements and 
alternative proofs, which potentially apply to more general situations 
(see Remark~\ref{r:propappl} below). We provide the alternative proofs in 
the following proposition, which is not used in the rest of the paper.

\begin{prop}\label{p:stab2}
Let $\ka_-$, respectively $\ka_+$, denote the negative, respectively
positive, Legendrian stabilization of the oriented Legendrian knot
$\ka\subset (Y,\xi)$. Then, for each $n>0$ we have:
\begin{enumerate}
\item
$\xi_{n+1}^-(\ka_-)$ is isotopic to $\xi_n^-(\ka)$;
\item
$\xi_n^-(\ka_+)$ is overtwisted.
\end{enumerate}
\end{prop}

\begin{proof}
This simple proposition can be deduced using the foundational results of Ko Honda from~\cite{Ho}. 
We refer the reader to~\cite{Ho} for the necessary background in what follows.
Let us quickly go over the details of the contact surgery construction. 
The contact framing together with the orientation on $\ka$ determine an oriented basis $\mu,\la$ 
of the first integral homology group of the boundary of a standard neighborhood $\nu(\ka)$ 
of $\ka$. The basis determines identifications 
\[
\del(\nu(\ka))\cong \R^2/\Z^2,\quad -\del(S^3\setminus\nu(\ka))\cong\R^2/\Z^2.
\] 
The surgery is determined by a gluing prescribed, with repect to the above
identifications, by the matrix
\[
A=
\begin{pmatrix}
n & -1\\
1 & 0\\
\end{pmatrix}.
\]
The pull--back of the dividing set is determined by 
\[
A^{-1}
\begin{pmatrix}
0\\1
\end{pmatrix}
=
\begin{pmatrix}
1\\n
\end{pmatrix},
\]
so it has slope $n$ on $\del(\nu(\ka))$. Applying a diffeomorphism of
the solid torus $\nu(\ka)$ this slope can be changed to $n/(1+nh)$ for
any $h\in\Z$. Therefore we can normalized it to lie between $-1$ and
$-\infty$, obtaining slope $-n/(n-1)$. By~\cite{Ho} there are exactly
two choices of tight contact structures on the solid torus with this
boundary slope, corresponding to the two possibile choices (positive
or negative) of a basic slice with boundary slopes $-1$ and
$-n/(n-1)$. With our conventions, choosing the negative basic slice
gives rise to the contact structure $\xi_n^-$. The knots $\ka_+$ and
$\ka_-$ can both be realized inside the neighborhood $\nu(\ka)$. If
$\nu(\ka_\pm)\subset\nu(\ka)$ is a standard neighborhood of $\ka_\pm$,
the (closure of the) difference $\nu(\ka)\setminus\nu(\ka_\pm)$ is a
basic slice, which is positive for $\ka_+$ and negative for
$\ka_-$~\cite{Ho}. Moreover, its boundary slopes with respect to the
basis $\mu, \la$ are $-1$ on $\del\nu(\ka_\pm)$ and $\infty$ on
$\del\nu(\ka)$.

For each $n\geq 0$, we can perform contact $(n+1)$--surgery along $\ka_\pm$ viewed as a Legendrian 
knot inside $\nu(\ka)$, obtaining another contact solid torus $\T$ with convex boundary in standard form. 
$H_1(\del\nu(\ka_-);\Z)$ has a basis $\mu'$, $\la'$ such that, with the obvious identifications, $\mu'=\mu$ 
and $\la'=\la-\mu$. Thus, since $\la=\mu'+\la'$, the identity 
\[
\begin{pmatrix}
n+1 & -1\\
1 & 0
\end{pmatrix}^{-1}
\begin{pmatrix}
0 & 1\\
1 & 1
\end{pmatrix}
=
\begin{pmatrix}
1 & 1\\
n+1 & n
\end{pmatrix}
\]
implies that, up to applying a diffeomorphism of $\T$, the slopes of $\del\nu(\ka)$ and $\del\nu(\ka_\pm)$ 
can be assumed to be, respectively,  $-n/(n-1)$ and $-(n+1)/n$. This shows that $\T$ can be 
decomposed as 
\[
\T = N\cup B,
\]
where $N$ is standard neighborhood of a Legendrian curve with slope
$-1$, and $B\cong T^2\x [0,1]$ has boundary slopes $-1$ and $-n/(n-1)$
and can be written as a union of two basic slices $B=B_1\cup B_2$,
where $B_1$ has boundary slopes $(-1,-(n+1)/n)$ and
$B_2=\nu(\ka)\setminus\nu(\ka_\pm)$, with the boundary slopes given
above. Since $B_2$ is a basic slice, $B$ is a basic slice (i.e.~it is
tight) if and only if $B_1$ and $B_2$ have the same sign as basic
classes~\cite{Ho}. By definition, $\xi^-_{n+1}(\ka_\pm)$ is the
contact structure obtained by taking $B_1$ to be a~\emph{negative}
basic slice.  Since $B_2$ is positive for $\ka_+$ and negative for
$\ka_-$, the analysis above proves simultaneously (1) and (2) of the
statement.
\end{proof}

\begin{rem}\label{r:propappl}
While the proof of  Proposition~\ref{p:stab} 
only holds, as written, for closed contact 3--manifolds, both the statement and 
the proof of Proposition~\ref{p:stab2} can stay the same 
even if $(Y,\xi)$ is open or has non--empty boundary. This allows one, at least 
in principle, to apply the approach of this paper in situations which are 
more general than the ones considered here. We hope to return to this issue 
in a future paper. 
\end{rem}

\begin{rem}\label{r:posstab}
Let $\xi_n^+(\ka)$ be the contact structure corresponding to the
choice of the~\emph{positive} stabilization $\ka_+$ of $\ka$. Then, an
argument analogous to that of the proof of Proposition~\ref{p:stab2} shows
that $\xi_{n+1}^+(\ka_+)$ is isotopic to $\xi_n^+(\ka)$ and
$\xi_{n+1}^+(\ka_-)$ is overtwisted.  Of course, this also follows
from the fact that $(-\ka)_-=-\ka_+$ and Proposition~\ref{p:stab2}.  In
fact,
\[
\xi_{n+1}^+(\ka_+) = \xi_{n+1}^-(-\ka_+) = \xi_{n+1}^-((-\ka)_-) 
= \xi_n^-(-\ka) = \xi_n^+(\ka)
\]
and 
\[
\xi_{n+1}^+(\ka_-) = \xi_{n+1}^-(-\ka_-) = \xi_{n+1}^-(\ka_+)
\]
\end{rem}

\section{The invariants: definition and basic properties}\label{s:def}

In this section we define the invariants, we prove some of their
properties and present some examples.

\subsection{Definition of the geometric invariant $I(\xi , \ka , f)$}
\label{ss:cosu}
Let $K$ be the knot type of a Legendrian knot $\ka\subset (Y,\xi)$,
and let $t:=\tb(\ka)\in\FF_K$ be the Thurston--Bennequin invariant of
$\ka$, i.e.~the contact framing of $\ka$. Suppose that either $K$ is
null--homologous or $Y$ has no $S^1\x S^2$--summand. For each positive
integer $n$, the contact structure $\xi_n^-(\ka)$ lives on the closed
3--manifold $Y_{t+n}(K)$ obtained by performing topological surgery
along $K$ corresponding to the framing $t+n$.

Let $(\ka,f)$ be a framed, oriented Legendrian knot in the contact
3--manifold $(Y,\xi)$, and let $\ka'\subset (Y,\xi)$ be a Legendrian
knot obtained by negatively stabilizing $\ka$ sufficientely many
times, so that $\tb(\ka')<f$. In view of Proposition~\ref{p:stab}, the
isomorphism class of the contact structure $\xi_{f-\tb(\ka')}^-(\ka')$ on
$Y_f(K)$ does not depend on the choice of $\ka'$ as long as
$\tb(\ka')<f$, therefore we can introduce the following:

\begin{defn}\label{d:Ileginvariant}
Assume that either the knot type $K$ is null--homologous or $Y$ 
is not of the form $Y'\# S^1\x S^2$. Let $(\ka,f)$ be a framed, oriented Legendrian knot in the contact 
3--manifold $(Y,\xi)$ such that $\ka$ has topological type $K$.  Define $I(\xi,\ka,f)$ to be the isomorphism 
class of the contact structure $\xi_{f-\tb(k')}^-(\ka')$ on $Y_f(K)$, where $\ka'\subset (Y,\xi)$ 
is any negative stabilization of $\ka$ such that $\tb(\ka')<f$. 
\end{defn}

\begin{prop}\label{p:sufflarge} 
Assume that either $K$ is null--homologous or $Y$ 
is not of the form $Y'\# S^1\x S^2$. Let $(\ka,f)$ be a framed, oriented Legendrian knot in the contact 
3--manifold $(Y,\xi)$ such that $\ka$ has topological type $K$. Then, $I(\xi,\ka,f)$ is overtwisted for 
each $f\leq\tb(\ka)$.  
\end{prop}

\begin{proof} 
By definition, $I(\xi,\ka,f)$ is the isomorphism class of the contact structure 
$\xi^-_{f-\tb(\ka')}(\ka')$, where $\ka'$ is any negative stabilization of $\ka$ 
such that $\tb(\ka')<f$. If $f\leq\tb(\ka)$ we can choose $\ka'$ 
so that $\tb(\ka')=f-1$. We have $\ka'=\ka''_-$ for some oriented Legendrian 
knot $\ka''$. Then,  
\[
\xi^-_{f-\tb(\ka')}(\ka')=\xi_1(\ka''_-)=\xi_1(-(-\ka'')_+)
\]
is overtwisted by Lemma~\ref{l:overtwisted}.
\end{proof} 

Recall that transverse knots admit a preferred orientation and can be
approximated, uniquely up to negative stabilization, by oriented
Legendrian knots~\cite{EFM, EH}.  Fix a transverse knot $\ta\subset
(Y,\xi)$, and let $\ka$ be a Legendrian approximation of $\ta$. Then,
by~\cite{EFM, EH}, up to negative stabilizations the Legendrian knot
$\ka$ only depends on the transverse isotopy class of $\ta$.

It follows immediately from Proposition~\ref{p:stab} that if $\ka'\subset (Y,\xi)$ is 
a negative stabilization of the oriented Legendrian knot $\ka\subset (Y,\xi)$, 
then for each framing $f$ we have $I(\xi,\ka',f)=I(\xi,\ka,f)$. This observation 
allows us to give the following: 

\begin{defn}\label{d:Itransvinvariant}
  Assume that either $K$ is null--homologous or $Y$ is not of the form
  $Y'\# S^1\x S^2$. Let $(\tau,f)$ be a framed transverse knot in the
  contact 3--manifold $(Y,\xi)$ such that $\tau$ has topological type
  $K$.  Define $\tilde I(\xi,\tau,f):=I(\xi,\ka,f)$, where $\ka$ is
  any Legendrian approximation of $\tau$.
\end{defn}

\begin{proof}[Proof of Theorem~\ref{t:maintransverse}]
Since the choice of $\ka$ is unique up to negative stabilization, 
the repeated application of Proposition~\ref{p:stab} verifies the result.
\end{proof}

\subsection{Heegaard Floer invariants}
We now apply the Heegaard Floer contact invariant defined by Ozsv\'ath
and Szab\'o \cite{OSz-cont}.

\begin{defn}\label{d:cleginvariant}
Let $Y$ be a closed, oriented 3--manifold, $K\subset Y$ a knot type and $f\in\FF_K$.
Assume that either $K$ is null--homologous or $Y$ is not of 
the form $Y'\# S^1\x S^2$. 
Given an oriented Legendrian knot $\ka\subset (Y,\xi)$, define
\[ 
c(\xi,\ka,f):=c(I(\xi,\ka,f))\in \wHF(-Y_f(K)), 
\]
and given a transverse oriented knot $\tau\subset (Y,\xi)$, define 
\[
\tilde c(\xi,\tau,f):=c(\tilde I(\xi,\tau,f))\in \wHF(-Y_f(K)).
\]
\end{defn}

\begin{rems}
\begin{itemize}
\item
It follows immediately from the definition, Lemma~\ref{l:overtwisted} and Proposition~\ref{p:stab} 
that, for each $f\in\FF_K$, $c(\xi,\ka_-,f)=c(\xi,\ka,f)$ and $c(\xi,\ka_+,f)=0$.
\item
It follows from Proposition~\ref{p:sufflarge} that $c(\xi,\ka,f)=0$ for 
each $f\leq \tb(\ka)$.
\item
If the complement of a Legendrian knot $\ka$ in $(Y,\xi)$ is overtwisted 
or has positive Giroux torsion, the same holds for $\xi_n^-(\ka')$ 
for some stabilization $\ka'$ of $\ka$. Therefore, it follows from the 
results of~\cite{GHV, OSz-cont} that $c(\xi,\ka,f)=0$ for each $f\in\FF_K$. 
\end{itemize}
\end{rems}


The following examples show that the invariant $c(\xi,\ka,f)$ is non--trivial.

\begin{exa}\label{exa:egy}
{\rm Consider the Legendrian unknot $\ka\subset (S^3,\xi_{st})$ with
Thurston--Bennequin number $-1$.  (In this case $\ka=-\ka$, so we
do not need to specify the orientation). Since the result of contact $(+1)$--surgery is 
equal to the unique Stein fillable contact structure on $S^1\times S^2$, we get that
$c(\xi_{st},\ka,\tb(\ka)+1)\neq 0$ (cf.~\cite[Lemma~5]{LS0}).}
\end{exa}

\begin{exa}\label{exa:ketto}
{\rm Let $\ka\subset (S^3,\xi_{st})$ be an oriented Legendrian knot with
knot type $K$ such that
\begin{equation}\label{e:tb}
\tb(\ka) = f_S+2g_s(K)-1>0,
\end{equation}
where $f_S$ is the framing defined by a Seifert surface of
$K$, and $g_s(K)$ is the slice genus of $K$. Then, by~\cite[Proof of
  Theorem~1.1]{LS1} $c(\xi_{st},\ka,\tb(\ka)+1)\neq 0$. As remarked
in~\cite{LS1}, the knot types containing Legendrian knots which
satisfy Condition~\eqref{e:tb} include all non--trivial algebraic
knots, i.e.~non--trivial knots which are links of isolated curve
singularities, as well as negative twist knots.}
\end{exa}


\subsection{The inverse limit construction}
The invariants $\tilde c(\xi , \ta , f)$ can be conveniently organized
as a single element in the inverse limit of certain Heegaard Floer
homology groups. In the rest of this section we spell out the details
of this construction.

Let $Y$ be a closed, oriented 3--manifold and $K$ a knot type in $Y$.
To each framing $f\in\FF_K$ one can naturally associate a triangle of
3--manifolds and cobordisms (cf.~\cite[pp.~933--935]{LS1}). Let
$Y_{f-1}(K)$ be the 3--manifold resulting from surgery along $K$ with
framing $f-1$. The first manifold in the triangle is $Y$, the second
is $Y_{f-1}(K)$ and the third one is $Y_f(K)$. A cobordism $W_f$ from
$Y_{f-1}(K)$ to $Y_f (K)$ can be given by considering a normal circle
$N$ to $K$ in $Y$, equip if with framing $f_S-1$, and after the
surgery on $K$ with framing $f-1$ has been performed, attach a
4--dimensional 2--handle along $N$ with the chosen framing. Simple 
Kirby calculus shows that $W_f$ is indeed a cobordism between 
$Y_{f-1}(K)$ and $Y_f(K)$.

When viewed upside down, $W_f$ induces a map 
\[
\widehat F_{\overline W_f}\co \wHF(-Y_f(K))\to \wHF(-Y_{f-1}(K)).
\]
Given framings $f\geq g$, define $\varphi_{g,f}$ to be the identity on $\wHF(-Y_f(K))$ if 
$f=g$, and the composition 
\[
\wHF(-Y_f(K))\xrightarrow{F_{\overline W_f}}\wHF(-Y_{f-1}(K))\lra\cdots
\xrightarrow{F_{\overline W_{g+1}}}\wHF(-Y_g(K))
\]
if $f>g$. Then, it is easy to check that the set 
\[
\left\{\left(\wHF(-Y_f(K)), \varphi_{g,f}\right)\right\}
\]
is an inverse system of $\Z/2\Z$--vector spaces and linear maps over the set $\FF_K$, 
so we can form the inverse limit $\Z/2\Z$--vector space 
\[
H(Y,K):=
\varprojlim\wHF(-Y_f(K)),
\]
which is the subspace
\[
\{(x_f)\in\prod_{f\in\FF_K}\wHF(-Y_f(K))\ |\ x_g=\varphi_{g,f}(x_f)\
\text{for $g\leq f$}\} \subset\prod_{f\in\FF_K}\wHF(-Y_f(K)).
\]

We define $\tilde c(\xi,\ta)$ as the vector
$\left(c(\tilde I(\xi,\ta,f))\right)_{f\in\FF_K}$, which is, \emph{a priori}
an element of $\prod_{f\in\FF_K}\wHF(-Y_f(K))$.
\begin{prop}\label{p:invlimit}
  The invariant $\tilde c(\xi,\ta)$ is in $H(Y,K)$.
\end{prop}
\begin{proof}
  Choose $f\in\FF_K$ and a negative stabilization $\ka'$ of
  $\ka$, with $\tb(\ka')<f-1$. By~\cite{DG}, performing contact
  $(+1)$--surgery on an extra push--off copy of $\ka'_-$ in the
  contact surgery presentation for $\xi_{f-\tb(\ka')}^-(\ka')$ gives
  $\xi_{f-1-\tb(\ka')}^-(\ka')$. The corresponding 2--handle attachment
  gives an oriented 4--dimensional cobordism from $Y_f(K)$ to
  $Y_{f-1}(K)$, and it is easy to check that reversing the orientation
  of that cobordism gives exactly the oriented cobordism $\overline
  W_f$. By~\cite[Theorem~2.3]{OSz-cont} (see
  also~\cite[Theorem~2.2]{LS1}), we have
\[
\widehat F_{\overline W_f}(c(\xi_{f-\tb(\ka')}^-(\ka'))) = c(\xi_{f-1-\tb(\ka')}^-(\ka')),
\]
i.e.
\[
\varphi_{f-1,f}(c(I(\xi,\ka,f)))=c(I(\xi,\ka,f-1)).
\]
Since this holds for each $f\in\FF_K$, and for $g\leq f$ we have 
\[
\varphi_{g,f} = \varphi_{g+1,g}\circ\cdots\circ\varphi_{f-1,f},
\]
the statement is proved.
\end{proof}

Proposition~\ref{p:invlimit} immediately gives the following:

\begin{cor}\label{c:defin}
If $c(\xi,\ka,g)\neq 0$ then $c(\xi,\ka,f)\neq 0$ for every $f\geq g$.
\qed
\end{cor}

\section{Proof of Theorem~\ref{t:main}}\label{s:main}

Let $\Si$ be an oriented surface-with-boundary and $\phi\co\Si\to\Si$
an orientation--preserving diffeomorphism which restricts to the
identity on a collar around $\del\Si$. Let $(Y,\xi_{(\Si,\phi)})$ be a
contact 3--manifold compatible with the open book decomposition
induced by $(\Si,\phi)$. Let $\tau\subset (Y,\xi_{(\Si,\phi)})$ be a
component of the boundary of $\Si$ viewed as the binding of the open
book and let $f_\Si$ be the framing induced on $\tau$ by $\Si$.

\begin{prop}\label{p:cobident}
There exists a Legendrian approximation $\ka\subset (Y,\xi)$ to $\tau$
such that $\tb(\ka)=f_\Si-1$ and, for each $n>0$, the contact
structure $\xi^-_n(\ka)$ admits a compatible open book with a binding
component $\tau'$ having the following properties:
\begin{itemize}
\item
Capping--off $\tau'$ gives back the open book $(\Si,\phi)$;
\item
Let $Z$ be the cobordism corresponding to capping--off $\ta'$, and let 
$X_{\ka, n}$ be the topological cobordism obtained by attaching a 
4--dimensional 2--handle along $\ka$ with framing $\tb(\ka)+n$. 
Then, $\overline Z=-X_{\ka,n}$, i.e.~$Z$ is obtained from $X_{\ka, n}$ 
by viewing it upside--down and reversing its orientation.
\end{itemize}
\end{prop} 

\begin{proof} 
As shown in~\cite[Lemma~3.1]{SVV}, any open book decomposition can be Giroux
stabilized so that the page of the new open book $(\Si',\phi')$ 
contains a Legendrian approximation $\ka$ of $\tau$ as a curve sitting
on a page $\Si'$ and parallel to (a component of) $\del\Si'$, with
$\tb(\ka)=f_{\Si'}=f_\Si-1$. This is illustrated in
Figure~\ref{f:approx}a--b.
\begin{figure}[ht!]
\labellist
\small\hair 2pt
\pinlabel $\Si$ at 10 241
\pinlabel $\Si'$ at 607 233
\pinlabel $\Si''$ at 616 73
\pinlabel $\tau$ at 147 292
\pinlabel $\tau$ at 331 294
\pinlabel $\tau$ at 205 154 
\pinlabel $+$ at 540 260 
\pinlabel $+$ at 344 96
\pinlabel $+$ at 521 91
\pinlabel $\ka$ at 385 252
\pinlabel $\ka$ at 324 53
\pinlabel $\ka_-$ at 199 90 
\pinlabel $\text{\scriptsize (a)}$ at 83 175 
\pinlabel $\text{\scriptsize (b)}$ at 450 175
\pinlabel $\text{\scriptsize (c)}$ at 365 4
\endlabellist
\centering
\includegraphics[scale=0.5]{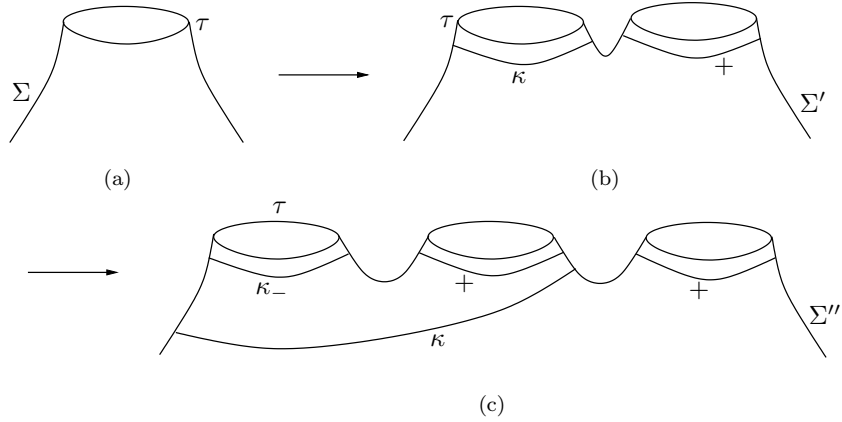}
\caption{Legendrian approximation of the binding. Once again, the
  signs $\pm $ on the curves indicate whether right-- or left--handed
  Dehn twists are to be performed on the given curve.}
\label{f:approx}
\end{figure}
As shown in Figure~\ref{f:approx}b--c,  after a further
Giroux stabilization both $\ka$ and its negative stabilization $\ka_-$
can be Legendrian realized on the page of an open book
$(\Si'',\phi'')$. The knot $\ka_-$ is parallel on $\Si''$ to a
boundary component which coincides with $\tau$, 
as indicated in Figure~\ref{f:approx}c.
Performing contact $(+n)$--surgery along $\ka$ is equivalent to a
contact $(+1)$--surgery along $\ka$ plus contact $(-1)$--surgeries
along $n-1$ parallel copies of $\ka_-$. Therefore, the resulting
contact structure is supported by the open book obtained by composing
$\phi''$ with a negative Dehn twist along the curve corresponding to
$\ka$, as well as positive Dehn twists along $n-1$ parallel copies of
the curve corresponding to $\ka_-$.  This is illustrated in
Figure~\ref{f:surg}a.
\begin{figure}[ht!]
\labellist
\small\hair 2pt
\pinlabel {\scriptsize $+$} at 93 257
\pinlabel {\scriptsize $+$} at 93 225
\pinlabel {\scriptsize $+$} at 253 240
\pinlabel {\scriptsize $+$} at 416 237
\pinlabel {\scriptsize $+$} at 323 82
\pinlabel {\scriptsize $-$} at 191 200
\pinlabel $\tau'$ at 72 279
\pinlabel $\text{\scriptsize n-1}$ at 170 269
\pinlabel $\Si'$ at 402 74
\pinlabel $\Si''$ at 492 245
\pinlabel $\text{\scriptsize (a)}$ at 254 170
\pinlabel $\text{\scriptsize (b)}$ at 254 13
\endlabellist
\centering
\includegraphics[scale=0.5]{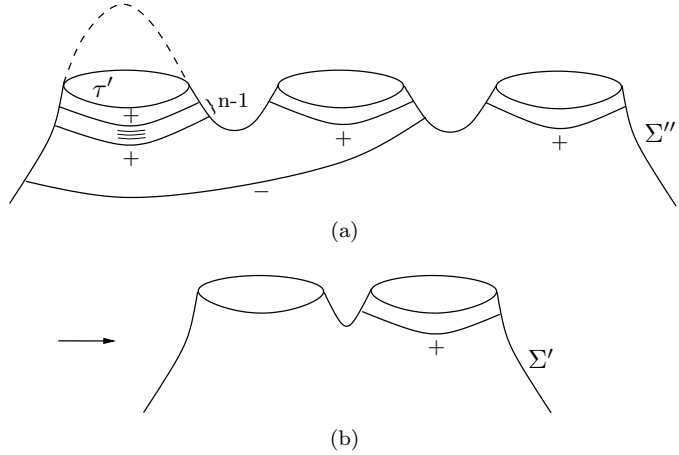}
\caption{Contact $n$--surgery and capping--off.}
\label{f:surg}
\end{figure}
As shown in Figure~\ref{f:surg}a--b, capping--off the
binding component denoted $\tau'$ in the picture yields the open book
$(\Si',\phi')$. (Notice the cancellation of the left--handed Dehn
twist of Figure~\ref{f:surg}a after the capping--off.)  This proves
the first part of the statement.

In order to control what happens at the level of 4--dimensional
2--handles, we represent the surgeries inside a standard Legendrian
neighborhood of $\ka$, as illustrated in Figure~\ref{f:neighbor}a.
\begin{figure}[ht!]
\labellist
\small\hair 2pt
\pinlabel {\Huge $\}$} at 306 354
\pinlabel $n-1$ at 340 355
\pinlabel {\Huge $\}$} at 298 105
\pinlabel $n-1$ at 332 106
\pinlabel $\tau'$ at 306 404 
\pinlabel $\tau'$ at  298 148
\pinlabel {\scriptsize $-2$} at  267 382
\pinlabel {\scriptsize $-2$} at  267 324
\pinlabel {\scriptsize $+1$} at  237 311
\pinlabel {\scriptsize $0$} at  264 155
\pinlabel {\scriptsize $-1$} at  264 133
\pinlabel {\scriptsize $-1$} at 264 75
\pinlabel {\scriptsize $+2$} at  215 62
\pinlabel {\scriptsize (a)} at 141 263
\pinlabel {\scriptsize (b)} at 141 15
\pinlabel $-1$ at 177 367
\endlabellist
\centering
\includegraphics[scale=0.55]{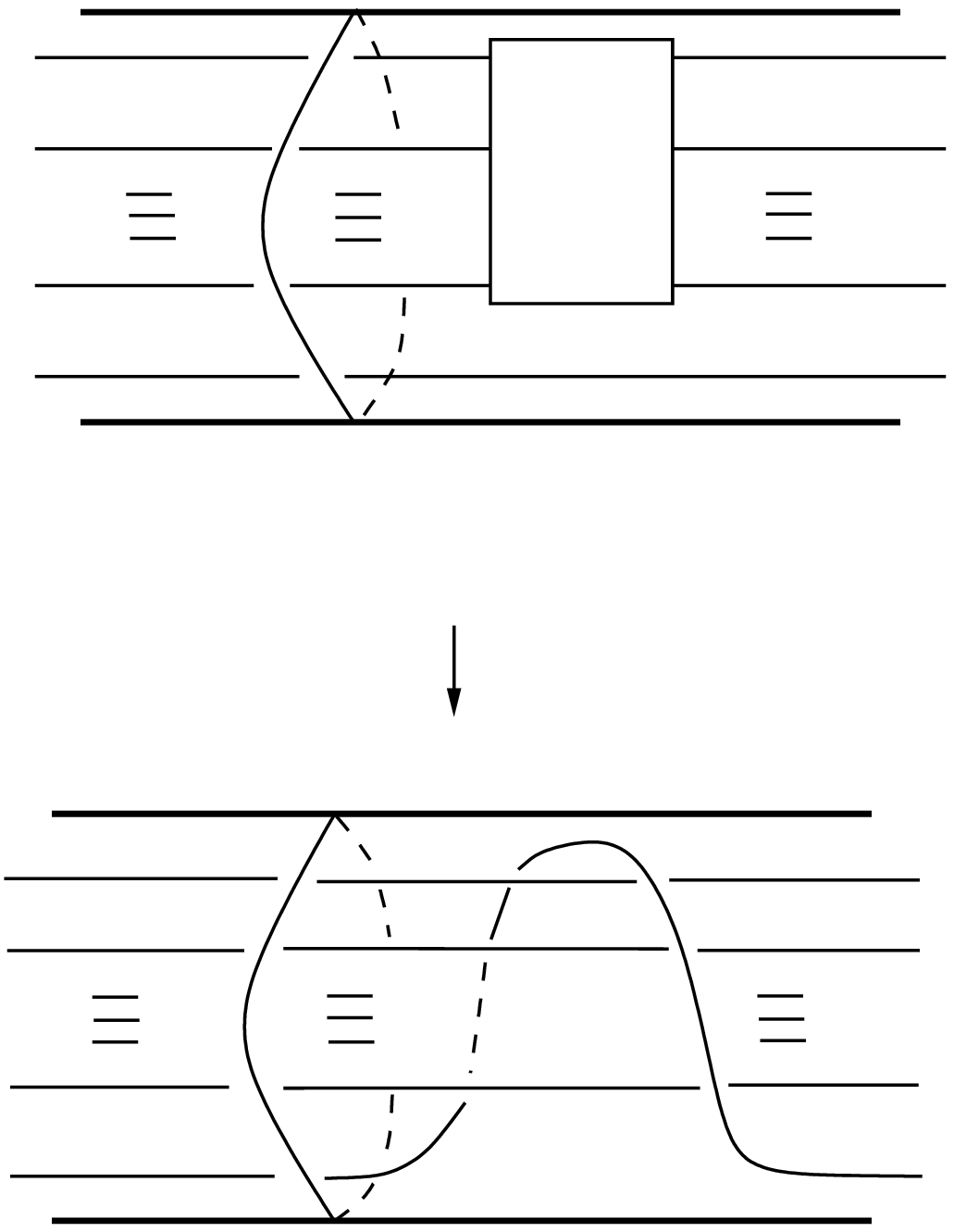}
\caption{Contact $(+n)$--surgery along $\ka$ in a standard neighborhood.}
\label{f:neighbor}
\end{figure}
\begin{figure}[ht!]
\labellist
\small\hair 2pt
\pinlabel $\tau'$ at  226 369
\pinlabel {\scriptsize $-2$} at  167 365 
\pinlabel {\scriptsize $-2$} at  153 321 
\pinlabel {\scriptsize $-2$} at  231 313
\pinlabel {\scriptsize $-2$} at  510 343
\pinlabel {\scriptsize $-2$} at  510 313
\pinlabel {\scriptsize $-1$} at  125 330
\pinlabel {\scriptsize $-1$} at  563 355
\pinlabel {\scriptsize $-1$} at  503 297
\pinlabel {\scriptsize $+1$} at  231 276
\pinlabel {\scriptsize $+1$} at  595 294
\pinlabel {\scriptsize $0$} at  139 91
\pinlabel {\scriptsize $0$} at  521 91
\pinlabel {\small $n$} at  214 71  
\pinlabel {\small $n+1$} at 588 73 
\pinlabel {\scriptsize (a)} at  125 236 
\pinlabel {\scriptsize (b)} at  530 236 
\pinlabel {\scriptsize (c)} at  530 13
\pinlabel {\scriptsize (d)} at 125 13
\endlabellist
\centering
\includegraphics[scale=0.58]{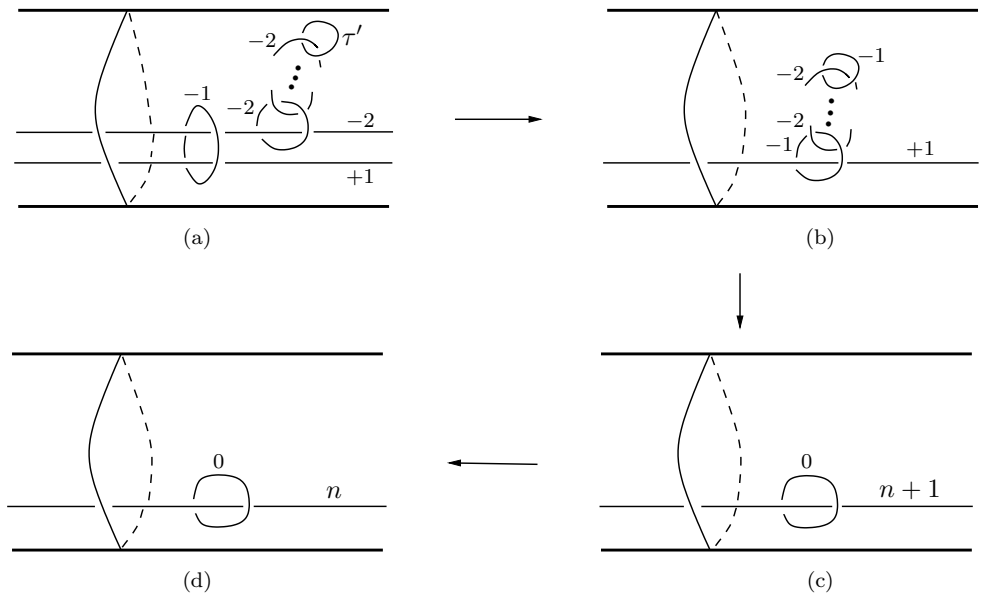}
\caption{Kirby moves in a standard neighborhood. In~\ref{f:kirby}a there are
  $n-2$ small $(-2)$--framed circles and a further long
  $(-2)$--framed circle. In~\ref{f:kirby}b we have $n-2$ $(-2)$--framed circles.}
\label{f:kirby}
\end{figure}
The framing coefficients appearing in the picture have the following
significance. Observe that each curve isotopic to the core of the
solid torus has a canonical framing coming from the identification of
the solid torus with $S^1\x D^2$. We have chosen the identification of
a neighborhood of $\ka$ with $S^1\x D^2$ so that the canonical framing
of the core corresponds to the contact framing of $\ka$. With this 
convention, the framing induced by $\Si''$ on $\tau'$ would be 
denoted by $``-1"$ in Figure~\ref{f:neighbor}a.
Figure~\ref{f:neighbor}b describes the surgery when 
we change the solid torus identification by a ``right--handed
twist''.  Then, the cobordism $Z$ corresponding to capping--off $\tau'$ is
obtained by attaching a 4--dimensional 2--handle along $\tau'$ with
framing $0$, as shown in Figure~\ref{f:neighbor}b. After sliding
$\ta '$ over one of the $(-1)$--framed circles in
Figure~\ref{f:neighbor}b, and then repeatedly the $(-1)$--framed
circles over each other, and finally blowing up the last two curves,
we arrive at Figure~\ref{f:kirby}a. A handle slide, followed by a
blow--down gives Figure~\ref{f:kirby}b, and further blow--downs give
Figure~\ref{f:kirby}c. Applying a ``left--handed twist'' to the
solid torus neighborhood (which just undoes the ``right--handed twist"
we applied earlier) gives Figure~\ref{f:kirby}d. This shows that $Z$
can be viewed as the cobordism obtained by attaching a 2--handle along
a meridian to the original curve $\ka$, with framing $0$ with respect
to the meridian disk. The picture shows that $Z$ coincides precisely
with $-\overline{X_{\ka,n}}$, where $X_{\ka,n}$ is the cobordism
obtained by attaching a 2--handle along $\ka$ with framing $+n$ with
respect to the contact framing, the minus sign denotes
orientation--reversal and the overline bar means viewing the cobordism
``up--side down''.
\end{proof}

\begin{proof}[Proof of Theorem~\ref{t:main}] 
By Baldwin's theorem~\cite[Theorem~1.2]{Ba}, for each positive integer
$n$ there is a Spin$^c$ structure $\s_0$ on the cobordism $\overline
Z$ such that
\[
F_{{\overline Z},\s_0}(c(\xi)) = c(\xi^-_n(\ka)).
\]

This equation then clearly shows that if $c(\xi)=0$ then $\tilde
c(\xi,\tau,f)=0$ for each $f\geq f_\Si$, verifying that the element
$\tilde c(\xi,\ta)\in H(Y, K)$ has only vanishing components, hence
$\tilde c(\xi,\ta)=0$.
\end{proof}

\section{Proofs of Theorem~\ref{t:mainsl} and Theorem~\ref{t:main2}}\label{s:main2}

Let $Y$ be a closed, oriented rational homology 3--sphere, and let
$\xi$ be a contact structure on $Y$. Let $\ka\subset (Y,\xi)$ be an
oriented Legendrian knot.  Suppose that $\ka=\del\Si$, where
$\Si\subset Y$ is an embedded oriented surface.  

Let $X$ be the oriented 4--dimensional cobordism obtained
by attaching a 4--dimensional 2--handle $H$ to $Y$ along $\ka$. 
Let $[\Si\cup D]\in H_2(X;\Z)$ denote the homology class 
supported by the union of $\Si$ and the core $D$ of
$H$.

\begin{lem}\label{l:spincunique}
If $\s_1$ and $\s_2$ are Spin$^c$ structures on $X_{\ka,n}$ having the same 
restriction to $Y$ and satisfying
\[
\langle c_1(\s_1),[\Si\cup D]\rangle = \langle c_1(\s_2),[\Si\cup D]\rangle,  
\]
then $\s_1=\s_2$. 
\end{lem} 

\begin{proof} 
  The set of Spin$^c$ structures on $X$ which restrict to $Y$ as a fixed 
  Spin$^c$ structure is an affine space on
  $H^2(X,Y)$ (with integral coefficients). Since $Y$ is a rational
  homology sphere and $X$ is obtained up to homotopy by attaching a
  2--disk to $Y$, by excision we have $H_2(X,Y)\cong\Z$ and
  $H^2(X,Y)\cong\Z$.  The exact homology sequence of the pair $(X,Y)$
  shows that the map $i_*\co H_2(X)\to H_2(X,Y)$ is injective, and
  therefore $H_2(X)\cong\Z$, with generator $[\Si\cup D]$. The exact
  cohomology sequence shows that the restriction map $i^*\co
  H^2(X,Y)\to H^2(X)$ is injective, and the free part of $H^2(X)$ has
  rank $1$. The evaluation map 
\[
H^2(X)\to Hom(H_2(X);\Z)\cong\Z,\quad 
  \be\mapsto\langle\be,[\Si\cup D]\rangle
\]
 is surjective. Therefore,
  the composition of $i^*$ with the evaluation map is injective. If
  $\s_1$, $\s_2$ are two Spin$^c$ structures on $X$ with coinciding
  restriction to $Y$, then $\s_1-\s_2=\al$, where $\al$ belongs to the
  image of $i^*$, and $c_1(\s_1)-c(\s_2) = 2\al$. Therefore, if the
  evaluation map takes the same value on $c_1(\s_1)$ and $c_2(\s_2)$,
  it follows that $\al=0$, hence $\s_1=\s_2$.
\end{proof} 

Let tb denote the Thurston--Bennequin number of $\ka$ with respect to
$\Si$, and let rot be the rotation number of $\ka$ with respect to
$\Si$.  Fix $n>0$, and let $X_{\ka,n}$ be the oriented 4--dimensional
cobordism obtained by attaching a 4--dimensional 2--handle to $Y$
along $\ka$ with framing tb$+n$.

\begin{prop}\label{p:spincvalue}
There exists a Spin$^c$ structure $\s$ on $X_{\ka,n}$ such that:
\begin{enumerate}
\item 
$\s$ extends the Spin$^c$ structures induced on $\del X_{\ka,n}$ by
  $\xi_{st}$ and $\xi^-_n(\ka)$;
\item
$\frac14(c_1(\s)^2-3\si(X_{\ka,n})-2\chi(X_{\ka,n}))+1 =
  d_3(\xi^-_n(\ka)) - d_3(\xi_{st})$;
\item
$\langle c_1(\s),[\Si\cup D]\rangle ={\mbox {rot}}+n-1$.
\end{enumerate}
\end{prop}

\begin{proof}
We can view Figure~\ref{f:neighbor}a (ignoring the knot $\tau$) as
$S^3$ union $n$ 4--dimensional 2--handles. The sequence of
Figures~\ref{f:kirby}a,~\ref{f:kirby}b and~\ref{f:kirby}c shows
that in fact Figure~\ref{f:neighbor}a represents $\widehat
X_{\ka,n}:=X_{\ka,n}\#(n-1)\overline{\CP}^2$. By e.g.~\cite[Section~3]{DGS},
$\widehat X_{\ka,n}\#\CP^2$ carries an almost complex structure $J$
inducing 2--plane fields homotopic to $\xi_{st}$ and $\xi^-_n(\ka)$ on
its boundary.  We define $\s_J$ to be the associated $\Spin^c$
structure, and $\s:={\s_J}|_{X_{\ka,n}}$.  By construction, $\s$
extends the Spin$^c$ structures induced on $\del X_{\ka,n}$ by
$\xi_{st}$ and $\xi^-_n(\ka)$. This proves Part (1) of the
statement. By~\cite{DGS} we have
\begin{equation}\label{e:d3formula}
\frac14(c_1(\s_J|_{\widehat X_{\ka,n}})^2-3\si(\widehat
X_{\ka,n})-2\chi(\widehat X_{\ka,n})) + 1 = d_3(\xi^-_n(\ka)) -
d_3(\xi_{st}).
\end{equation}
Figure~\ref{f:neighbor}b gives a natural basis
$(\be,x_1,\ldots,x_{n-1})$ of $H_2(X_{\ka,n};\Z)$ satisfying $\be\cdot
x_1=1$ and $x_i\cdot x_{i+1}=1$ for $i=1,\ldots,n-2$. By construction
and~\cite{DGS}, the values of $c_1(\s_J)$ on this basis are given by
$\langle c_1(\s_J),\be\rangle =$rot, $\langle c_1(\s_J), x_i\rangle =$
rot$-1$, $i=1,\ldots, n-1$. We want to express the generator $[\Si\cup
  D]$ of $H_2(X_{\ka,n};\Z)\cong\Z$ in terms of $\be$ and $x_1,\ldots,
x_{n-1}$. In Figure~\ref{f:kirby}b, the classes represented by the
framed circles (except the $-1$--framed one corresponding to $\ta$)
give us the new basis
$(\be,\be-x_1,x_1-x_2,\ldots,x_{n-2}-x_{n-1})$. If we define the
classes $e_1,\ldots, e_{n-1}$ by setting
\[
e_1:=\be-x_1,\ e_{i+1}-e_i:=x_i-x_{i+1},\ i=1,\ldots, n-2,
\]
it is easy to check that 
\[
[\Si\cup D]=\be+e_1+\ldots,e_{n-1}
\]
and $\langle c_1(\s_J),e_i\rangle = 1$ for $i=1,\ldots, n-1$. Thus, 
\[
\langle c_1(\s), [\Si\cup D]\rangle = \langle c_1(\s_J), \Si\cup
D]\rangle = {\mbox {rot}}+n-1.
\]
This proves Part (3) of the statement. Finally, the values $\langle c_1(\s_J),e_i\rangle = 1$
imply that 
\[
c_1(\s_J|_{\widehat X_{\ka,n}})^2-3\si(\widehat X_{\ka,n})-2\chi(\widehat X_{\ka,n}) =
c_1(\s)^2-3\si(X_{\ka,n})-2\chi(X_{\ka,n}).
\]
Thus, Equation~\eqref{e:d3formula} implies Part (2) of the statement. 
\end{proof}

\begin{proof}[Proof of Theorem~\ref{t:main2}]
Let $Z$ be the cobordism of Proposition~\ref{p:cobident} corresponding 
to an integer $n$. By~\cite[Theorem~1.2]{Ba}, there 
is a Spin$^c$ structure $\s_0$ on the cobordism ${\overline {Z}}$ 
such that
\begin{equation}\label{e:bald}
F_{{\overline Z},\s_0}(c(\xi_{st})) = c(\xi^-_n(\ka)),
\end{equation}
where $\ka$ is the Legendrian approximation of $\tau$ described in 
Proposition~\ref{p:cobident}. Since $c(\xi_{st})\neq 0$, we will prove that 
$\tilde c(\xi_{st},\tau,f_S+2g(K))\neq 0$ by showing 
that for an appropriate choice of $n$ the map 
$F_{{\overline Z},\s_0}$ is injective. 

By Proposition~\ref{p:cobident}, the cobordism ${\overline Z}$ is 
exactly the cobordism induced by $-p$--surgery along the mirror 
image knot $\overline{K}$, where $p=\tb(\ka)+n$. 
Now we choose $n$ so that $p=2g(K)$. Then, the assumption 
$\sel(\tau)=\tb(\ka)-\rot(\ka)=2g(K)-1$ implies 
\[
\rot(\ka)+n-1=\tb(\ka)-2g(K)+1+2g(K)-\tb(\ka)-1=0.
\] 
Therefore, the $\Spin^c$ structure $\s$ of 
Proposition~\ref{p:spincvalue} satisfies $c_1(\s)=0$. 

By Equation~\eqref{e:bald} and the identification $\overline{Z}=-X_{\ka,n}$, 
the $\Spin^c$ structure $\s_0$ satisfies
\begin{equation}\label{e:degshift}
\frac14(c_1(\s_0)^2-3\si(-X_{\ka,n})-2\chi(-X_{\ka,n}))=-d_3(\xi^-_n(\ka))+d_3(\xi_{st}) 
\end{equation}
Since $\si(-X_{\ka,n})=-\si(X_{\ka,n})$, Equation~\eqref{e:degshift} together with 
Proposition~\ref{p:spincvalue}(2) imply 
$c_1(\s_0)^2=-c_1(\s)^2=0$, therefore $c_1(\s_0)=0$. 

Let $\t_0$ denote the restriction of $\s_0$ to $S^3_{-2g(K)}(\overline{K})$. 
By~\cite[Theorem~9.19 and Remark 9.20]{OSzF2}, there is a surjective map 
\[
Q\colon \Spin ^c (S_0^3(\overline{K}))\to \Spin^c (S_{-2g(K)}^3(\overline{K}))
\]
and an exact triangle 
\begin{equation}\label{e:triangle}
\begin{graph}(6,2)
\graphlinecolour{1}\grapharrowtype{2}
\textnode {A}(1,1.5){$\wHF (S^3)$}
\textnode {B}(5, 1.5){$\wHF (S^3_{-2g(K)}(\overline{K}),\t_0)$}
\textnode {C}(3, 0){$\wHF (S^3_0(\overline{K}), [\t_0])$}
\diredge {A}{B}[\graphlinecolour{0}]
\diredge {B}{C}[\graphlinecolour{0}]
\diredge {C}{A}[\graphlinecolour{0}]
\freetext (2.4,1.8){$F$}
\end{graph}
\end{equation}
where 
\[
\wHF (S_0^3(\overline{K}), [\t_0]):=
\oplus _{\t\in Q^{-1}(\t_0)}\wHF (S_0^3(\overline{K}),\t).
\]
We claim that 
for each $\t\in Q^{-1}(\t_0)\subset \Spin ^c (S_0^3({\overline {K}}))$ we have 
\[
\vert \langle c_1 (\t), h\rangle \vert \geq 2g(K),
\]
where $h$ is a homology class generating $H_2(S^3
_0(\overline{K});\Z)$. In fact, since $\s_0$ extends $\t_0$,
by~\cite[Lemma~7.10]{abs},
\[
0=\langle c_1(\s_0),[\Si\cup D]\rangle \equiv 
-2g(K) + \langle c_1(\t),h\rangle \qquad \bmod 4g(K),
\]
which immediately implies the claim. 

Notice that $h$ can be represented by a genus-$g(K)$ surface (by
adding the core of the 2--handle to a Seifert surface of
$\overline{K}$), therefore the adjunction formula of
\cite[Theorem~7.1]{OSzF2} implies $\wHF(S_0^3 ({\overline
  {K}}),[\t_0])=\{0\}$.  This shows that the horizontal map $F$
of~\eqref{e:triangle}, is an isomorphism.  Since $c(\xi _{st})$
generates $\wHF (S^3)$, we have $F(c(\xi_{st}))\neq 0$.  In view of
Equation~\eqref{e:bald}, to prove that $\tilde
c(\xi_{st},\tau,f_S+2g(K))\neq 0$ it suffices to show that
$F=F_{\overline{Z},\s_0}$. From the general theory we know that
\[
F=\sum _{\{\s\in Spin ^c(\overline{Z})\ \vert\  
\s|_{S^3 _{-2g(K)}(K)}=\t_0\}} F_{\overline{Z},\s}. 
\]
Since $F$ is an isomorphism on $\wHF (S^3)=\Z /2\Z$, all the $\Spin^c$
structures contributing nontrivially the sum have the same degree
shift. Moreover, if $\overline{\s}$ denotes the $\Spin^c$ structure
conjugate to $\s$, by~\cite[Theorem~3.6]{OSz-4d} we have
\[
F_{\overline{Z},\overline{\s}}(c(\xi_{st})) = J F_{\overline{Z},\s}(c(\xi_{st})), 
\]
where $J$ is the identification between $\wHF (Y, \s)$
and $\wHF (Y, {\overline {\s}})$ defined in~\cite{OSzF1}. 
Since $J$ preserves the absolute
grading, this implies that each $\Spin^c$ structure $\s$ with
$\s\neq\overline{\s}$ contributes trivially to
$F(c(\xi_{st}))$. Finally, from $c_1(\s_0)=0$ we know that
$\s_0=\overline{\s_0}$, therefore by Lemma~\ref{l:spincunique} $\s_0$
is the only $\Spin$ structure on $\overline{Z}$ which extends
$\t_0$. This implies $F=F_{\overline{Z},\s_0}$ and proves $\tilde
c(\xi_{st},\tau,f_S+2g(K))\neq 0$. Applying Corollary~\ref{c:defin} it
follows that $\tilde c(\xi_{st},\tau,f)\neq 0$ for each $f\geq
f_S+2g(K)$.
\end{proof}

\begin{proof}[Proof of Theorem~\ref{t:mainsl}]
Let $\tau$ be a transverse knot in the knot type $K$ with $\sel (\ta
)=2g(K)-1$.  According to~\cite[Lemma~6.5]{BEVM} there is an open book
decomposition of $S^3$ compatible with $\xi _{st}$ having one binding
component equal to $\ta$. Applying Theorem~\ref{t:main2} the result
follows for each integer surgery coefficient $n\geq 2g(K)$. In
particular, Theorem~\ref{t:main2} gives a contact structure on
$S^3_{2g(K)}(K)$ with nonvanishing contact invariant.  By the
algorithm for contact surgeries described in Section~\ref{sec:second},
for $r\in\Q$ with $r\geq 2g(K)$ a contact structure can be given by
performing an appropriate sequence of Legendrian surgeries on the
contact structure previously constructed on $S^3_{2g(K)}(K)$. Since
under Legendrian surgery the nonvanishing property of the contact
invariant is preserved, the 3--manifolds $S^3_r (K)$ with $r\geq
2g(K)$ all carry contact structures with nonvanishing contact
Ozsv\'ath-Szab\'o invariants. Since $c(Y,\xi )\neq 0$ implies
tightness for $(Y, \xi)$, the proof of the theorem is complete.
\end{proof}

\end{document}